\input phyzzx\def\e{\adveq\eqno{\rm (\chapterlabel\the\equanumber)}}

\def\adveq{\global\advance\equanumber by 1}
\def\myeq{{\rm \chapterlabel\the\equanumber}}
\def\rarrow{\rightarrow}

\def\semidirect{\mathrel{\raise0.04cm\hbox{${\scriptscriptstyle |\!}$
\hskip-0.175cm}\times}}


\def\ref#1{$^{[#1]}$}

\def\r#1{$[\rm#1]$}

\def\e{\adveq\eqno{\rm (\chapterlabel\the\equanumber)}}

\def\adveq{\global\advance\equanumber by 1}
\def\myeq{{\rm \chapterlabel\the\equanumber}}
\def\rarrow{\rightarrow}

\def\semidirect{\mathrel{\raise0.04cm\hbox{${\scriptscriptstyle |\!}$
\hskip-0.175cm}\times}}


\def\ref#1{$^{[#1]}$}

\def\r#1{$[\rm#1]$}

%


\def\Conjecture{\par\noindent{\bf Conjecture\ }}
\def\WIS{\address{Department of Physics\break
Weizmann Institute\break Rehovot 76100, Israel\break}}
\overfullrule=0pt
\line{\hfill WIS}
\date{April, 2003}
\titlepage
\title{On the Distribution of Prime Multiplets}
\author{Doron Gepner}
\WIS
\abstract
The probability of finding a prime multiplet, i.e., a sequence of primes $p$ and $p+a_i$, 
$i=1\ldots m$,
being all primes where $p$ is some prime less than the integer $n$ is naively $1/log(n)^{m+1}$.
It is shown that, in reality, it is proportional to this probability by a constant factor which depends on $a_i$ and $m$ but not on $n$, for large $n$. These constants are
appellated as PDF (prime distribution factors). Moreover, it is argued that the PDF depend on the
$a_i$ in a "week" way, only on the prime factors of the differences $a_i-a_j$ and not
on their exponents. For example $p$ and $p+2^s$ will have the exact same probability for all
integer $s>0$. The exact formulae for
the PDF ratios are given. Moreover, the actual 'basic' PDF's are calculated exactly
and are shown to be less than $1$, which indicates that primes 'attract' each other.
An exact asymptotic formula for the number of basic multiplets is given. 
\endpage
\par
Prime numbers have arose curiosity for a long time. Their distribution received 
attention by the conjecture by Gauss that between $2$ and $n$ there are
on the average $n/\log(n)$ primes. This was proved in 1896 by Hadamard and de La Vall\'ee Poussin.  Therefore, the probability of finding a prime less than $n$
is $1/\log(n)$, for large enough $n$. Recently attention was given to multiplets of primes,
by the work of Goldstone and Yildrim \REF\GY{Goldstone and Yildrim, Nature.}\r\GY. For example, twin primes, such as $p$ and $p+2$ both
being primes. Or, sequences such as $p$, $p+2$, $p+6$ and $p+8$, all being prime.

In this letter I wish to make several observations on prime multiplets, described below.
Suppose that $a_i$, $i=1,..,m$ is a list of m positive integers, in ascending order. Then we may ask what is
the probability that $p$ and $p+a_i$ will all be primes, where $p$ is a prime less than $n$. Let us
denote the number of such primes by $N(n,a_i)$.
Statistically, if this would be independent events, the average number should be 
$$n/\log(n)^{m+1}.\e$$

The probability is, actually, different. Moreover the ratio between the measured probability
and the naive one is a constant, typical for the list of numbers $a_i$. This constant, 
$f(a_i)$,
is defined by
$$f(a_i)=n/N(n,a_i)/log(n)^{m+1},\e$$ 
where $n$ is large enough. If $N(n,a_i)$ is zero, we define $f(a_i)$
to be infinite.
For example, for $p$ and $p+2$ we find that the constant is $0.66$ 
approximately, with $n$ about $10^6$.
Another example is $p$, $p+2$, $p+6$. Here we find for $10^6$ that the factor
is approximately $0.27$.
We call $f(a_i)$ prime distribution factors or PDF. As can be seen later, these numbers
are off the limiting values, which are calculated later in this note, eq. (14). For example, the actual
limit for $p$ and $p+2$ is $0.757392...$, and not $0.66$.  To calculate the exact numbers, 
directly, one needs
to go to extremely large numbers. However, the ratios of PDF's for the same $m$ 
can be reliably estimated for a few million primes.

Indeed, for many different
$a_i$ that we tried, the ratio eq. (2) exist and can be measured easily.
This leads us to conjecture number (1):

\Conjecture{(1): }
 For all sequence $a_i$ the limit on the r.h.s. of eq. (2) exists and enables 
us to define the PDF, $f(a_i)$.

One may wander how does the PDF $f(a_i)$ depends on $a_i$. For example, $p$ and $p+a$ being 
both primes where $a$ is any even number. For this case a remarkable phenomena happens.
The PDF is dependent only on the list of prime factors in the decomposition of $a$. I.e.,
if 
$$a=\prod p_i^{s_i},\e$$ 
where
$p_i$ are some primes, then the PDF depends only on $p_i$ and not the exponents $s_i$.
For example, for $a=2,4,8$, etc, it is the same PDF, or for all $a=2^n$.
Similarly, for $a=2^b 3^c$ it is the same for any integers $b,c>0$.

One may wander then how this phenomena extends to larger chains, $m>1$.
It turns out to obey a simple rule too. The PDF $f(a_i)$ depends only on the primes
composing the number
$$x=\gcd_{i,j} \{a_i-a_j\},\e$$
where $i,j$ range on all possible differences $i,j=0\ldots m$ where $a_0=0$ is set by
convention, and on $m$. For example, $p$, $p+2$ and $p+6$ will have the same PDF as 
$p$, $p+4$ and $p+96$  since the gcd eq. (4) is composed only
from the primes $2$ and $3$, in both cases. Another example is $p$, $p+8$ and $p+12$, being
all primes. 
It is easy to enumerate many examples
of this kind. This leads us to the following second conjecture,

\Conjecture{(2): }
The PDF eq. (2), $f(a_i)$, depend only on the primes composing the gcd, eq. (4), where 
$a_0$ is set to zero, and $i,j=0\ldots m$.
 
To summarize, this is an interesting observation on the distribution of prime 
multiplets. One could expect to see the naive probability when the pairs are widely
apart, i.e., say, $p$ and $p+a$ where $a$ is large, and thus to become independent.
This is not at all the case. For $a=2^n$, according to conjecture $(2)$ it is always the
same probability. $f$ depends olny on the prime factors. This is 
definitely an indication that the distribution of primes, rather than being random,
is "non-locally" strongly correlated. More on the distribution later in this note.

Let us turn now to some examples. We start by considering twins, $m=1$,
of the form $p$ and $p+a$. These are some sample calculations.

We checked all the primes up to the $400000$th prime which is $5800079$.
For the case $a=2$ we find that there are exactly 36826 pairs, starting
with $3,5$ and $5,7$, etc. The ratio eq. (2) is then $f\approx 0.6494$.

For $a=4$ we find almost the same number, according to conjecture (2),
i.e., $36707$ pairs. For $a=6$ we find $73187$ pairs with $f\approx 0.32676$. For $a=12$, again it is almost the same number $73449$,
again verifying conjecture (2). For $a=14$ there are $43993$ pairs and $f=0.543606$. For $a=30$ we find 97825 pairs and $f=0.244466$.

Actually, if one considers the ratios of these numbers they are very close to simple rational numbers, e.g. $f(6)/f(2)\approx 0.5$ this leads us to an exact formula for the binary pairs,
$$f(a)=f(2)\prod_i {p_i-2\over p_i-1},\e$$
where $i$ ranges over all prime factors of the number $a$, except $2$.
We verified this formula for many cases, e.g. $a=70, 210, 30, 14$, etc.,
and it is exact to few tenth of a percent, for the first million primes.

This leads us to the third conjecture.

\Conjecture{(3): }
The ratios of the numbers $f(a_i)$ for a fixed $m$ are given by simple
rational numbers. For $m=1$ the ratio is given by eq. (5).

Let us turn now to examples of triplets. The simplest one is $p$, 
$p+2$ and $p+6$. The number $x$ is always divisible by two and three.
We find for the triplets from one to $5800079$ the factor $f_0=0.278193$.
For the triplet $p$, $p+2$ and $p+12$ we find $f=0.182965$.
For $p$, $p+2$ and $p+14$ we find $f=0.222554$. For $p$ $p+6$, $p+70$ 
we find $0.14695$. Again we note that the ratios of these numbers are
rational and that $f$ is given by
$$f(a_i)=f_0\prod_{p_i} {p_i-3 \over p_i-2} ,\e$$
where $p_i$ ranges over all prime divisors of $x$, eq. (4),
excluding $2$ and $3$. 

From these two cases we can guess the general formula for any $m$-plet.
It is given by

$$f(a_i)=C(m) \prod_{p_i}  {p_i-m-1\over p_i-m-1+g_i},\e$$
where $p_i$ ranges over all prime factors of $x$, eq. (4), which are greater than $m+1$, and $g_i$ is the number of independent differences
$a_i-a_j$ divisible by $p_i$.

This formula  is the conjecture:

\Conjecture{(4): }
The PDF are given by eq. (7).

Now, let us turn our attention to the constants $C(m)$. These are the PDF for a multiplet divisible
only by primes less or equal to $m+1$, where $m+1$ is the length of the multiplets.
These we term, basic multiplets.
For example for $m=1$, $N(x,m)$ is the number of pairs $p$ and $p+2$, which are both prime, $p\leq x$.
Interestingly, there is a conjecture for this number 
\REF\Br{R. P. Brent, Math. Comp., 28 (1974).}\r\Br, and ref. therein;
for a review see \REF\En{Encyclopedic Dictionary of Mathematics, 
Distribution of prime numbers, 128D.}\r\En, which is
$$N(x,1)\approx k \int_2^x {du\over log(u)^2},\e$$
where 
$$k=2\prod_{p>2\atop p {\rm\  prime}} \{ 1-{1\over (p-1)^2} \}=1.32032...\e$$
This conjecture was verified in calculations, e.g. for $x=10^9$ there are
$3424506$ pairs, agreeing very well with this formula which gives $3425230$.

As we have found, this formula enjoys a generalization to all basic $m$ plets, for any $m$.
The number $N(x,m)$ is conjectured to be approximately,
$$k(m)=z(m)\prod_{p>m+1 \atop p {\rm\  prime}} \{ 1-{1\over (p-q+1)^{m+1}} \},\e $$

where $q$ is the highest prime less or equal to $m+1$, and
$$N(x,m)=k(m) \int_{m+1}^x {du\over log(u)^{m+1}},\e$$
and where $z(m)$ is an integer conjectured to be
$$z(m)=(m+1)(m-2)\ldots (m+1-3t),\e$$
where $t$ is the highest integer such that $m+1-3t>0$.

Let us give several examples to this formula. For the basic triplet we take $p$, $p+2$ and
$p+6$. Up to the $400000$ prime there are 
$5520$ such triplets. From the formula it comes to $5580$, up to $x=5800079$.
Up to the $10^6$ primes there are $12092$ basic triplets, and from the formula, eq. (11),
$12170$. Up to the $2\times 10^6$ there are $21953$ triplets, and from the formula
$22099$, up to the prime $32452843$. This concludes the evidence for triplets,
giving credence to the formula eq. (11), for the case of triplets.

Eq. (11) can be checked also for quadruplets, $m=3$. Since there is less statistics, we expect
the convergence to be worse. 
For the basic quadruplet we take $p$, $p+2$, $p+6$ and $p+8$.
Up to the $400,000$th prime, we find 591 such quadruplets,
up to the prime $5800079$. The formula gives $551.54$. For the first $10^6$ primes we find
1229 primes up to $15485863$. Eq. (11) gives $1115.5$. We checked also $2\times 10^6$.
Here we find $2052$ primes up to $32452843$. The formula gives $1923$ primes, or it is about 
$5$\% off. We believe that bigger primes will indeed converge to the asymptotic equation (11).

Next we check quintuplets. Here we take, $p$, $p+2$, $p+6$, $p+8$ and $p+12$.  
We get good agreement with the asymptotic formula, eq. (11). For the
first $4\times 10^5$ primes we have $109$ such quintuplet, where the formula gives
$103$, up to the prime $5800079$. For the $10^6$ numbers we have $205$ multiplets, where 
the formula gives $191.36$, up to $15485863$. We checked also the first $10^7$ primes. We find
$336$ quintuplets, where the formula gives $311.6$, up to the prime $179424673$. 

Lastly, we checked the asymptotic formula for sextuplets. Here we take for the basic multiplet
$p$, $p+2$, $p+6$, $p+8$, $p+12$ and $p+18$. Up to $5800079$ there are
$15$ such sextuplets, whereas the formula gives $16.09$. Up to $15485863$ there are
$20$ multiplets, whereas the formula gives $25.99$. For $86028121$ we get
$57$ multiplets, whereas the formula gives $68.61$. 

These results 
encourage us to believe that with further
calculations a rather exact correspondence could be seen, and that the formula eq. (11) is 
asymptotically exact.  

Now, we come to the question of determining the basic PDF, $C(m)$, eq. (7). The function 
eq. (11) has the limit
$$\lim_{x\rarrow\infty} log(x)^{m+1}/x\int_{m+1}^x {du\over log(u)^{m+1}}=1,\e$$
as is easy to see by performing the integral by parts, 
giving this up to negligible 
pieces. This implies that the PDF $C(m)$ is 
$$C(m)=1/k(m),\e$$
Interestingly, $C(m)$ is less than one, e.g., $C(1)=0.757392...$,
$C(2)=0.34997\ldots$,
implying that these multiplets are more frequent than what may be naively expected. 
This shows that,
in fact, the
primes "attract" each other.

This forms our last conjecture:
\Conjecture{(5): }
The number of basic $m$-plets up to the number $x$ is given by eq. (11). $1/k(m)$ is the basic 
PDF, equal to $C(m)$, eq. (10), and it is always less than one.

There is actually a probabilistic way to understand eq. (11). 
Consider the
pairs $p$ and $p+2$. The probability of either being prime is $2/p$. 
By the sieve 
method then the number of prime pairs less or equal to $x$ is approximated by
$$M={x\over2}  \prod_{p>2\atop p {\rm\  prime}}^{\sqrt x}\{1-2/p\}.\e$$
 
This generalizes trivially
to the higher multiplets, where the probability is $(m+1)/p$, thus replacing $2$ with $m+1$. 
$$M(m)=x Z(m)\prod_{p>m+1 \atop p {\rm\ prime}}^{\sqrt x} \left\{1-{m+1\over p} \right\},\e$$
where 
$$Z(m)=\prod_{u\leq m+1\atop u \ \rm prime} {1\over u}.\e$$ 
In passing,
we note another way of expressing $M(m)$,
$$b(m)=x/M(m)=Z(m)^{-1}  \sum_{t} {(m+1)^{l(t)}\over t},\e $$
where $t$ is any number whose prime factors are all primes less than $\sqrt x$, including
of course all numbers up to $x$, which is not divisible by the primes less or equal to 
$m+1$, and 
$$l(\prod p_r^{s_r})=\sum s_r,\e$$
where $p_r$ are primes bigger than $m+1$.

Now, take a non-basic multiplet, e.g., $p$ and $p+2 q$, 
$q$ prime.
Then for the prime $q$ the probability changes to 
$$(1-1/q)/(1-2/q)=(q-1)/(q-2),\e$$ 
since it is enough to
check only one prime.  
This inverse ratio becomes
$${q-2\over q-1},\e$$
which is precisely the PDF we found, eq. (5).
 For several primes indeed the probability
is a product of all such factors.

For $m$ greater than one, the probabilistic argument indeed 
gives eq. (7). To see this, suppose that a given prime $p$ 
divides {\bf several} independent $a_{ij}=a_i-a_j$, in the notation of eq. (4). 
Denote the maximal number of such
divisible differences by $g$. From the probabilistic
argument it follows that the probability ratio is 
$$K=(1-(m+1-g)/p)/(1-(m+1)/p)=(p-m-1+g)/(p-m-1),\e$$
since it is enough to check only $m+1-g$ numbers divisable by $p$,
instead of the basic $m+1$ numbers,
and the PDF is
$$C(m)/K,\e$$
exactly verifying eq. (7).
We checked this in examples, and indeed it works, e.g., for $p$, $p+10$ and $p+30$,
we find that the ratio of PDF is $0.5$ in accordance with eq. (22).

Thus, it appears that a probabilistic argument indeed explains the values of the PDF's.
This is an indication that these values are probabilistic,
and that this simple sieve argument gives the exact values.

Actually, it is not difficult to compute the sieve product $M(m)$, 
eq. (16).
We need two identities proven by Martens (1874) \r\En:
$$\sum_{p\leq x \atop p \ \rm prime} {1\over p} =\log\log (x) +B+
O({1\over \log(x)}),\e$$
where $B$ is a constant equal to $B=0.2616...$ and
$$\prod_{p\leq x \atop \ \rm prime} \left( 1-{1\over p} \right)=
{e^{-c}\over \log(x)} \left( 1+O({1\over\log(x)})\right),\e$$
where $c$ is Euler's constant $c=0.577215...$.

Now, consider the product, eq. (25). It is the sieve for single primes,
i.e., the number of primes up to $x$ is given by
$$r_0\,x\, \prod_{p\leq \sqrt x \atop p \ \rm prime} \left( 1-{1\over p}
\right )={x\over \log(x)}\left[1+O({1\over \log(x)})\right],\e$$
where we used eq. (25) and set 
$$r_0=\exp(c)/2,\e$$
to get the correct result. $r_0$ is less, and close to one,
$r_0=0.890536$. We conclude that the sieve sum needs to be corrected
by a factor. So we redefine the sieve product,
$$M(m)=Z(m)\, x\, r_m \prod_{q> m+1 \atop q \ \rm prime}^{\sqrt x} \left ({1-{m+1\over q}}\right),\e$$
where $r_m$ is a factor close to one, as yet to be determined.

We can now take the log of eq. (28). We find that
$$H(m)=\log\left[{M(m)\over x\,Z(m)}\right]=\sum_{q>m+1 \atop q \ \rm prime}^{\sqrt x}
\log\left({1-{m+1\over q}}\right).\e$$
Expanding the Log in series we obtain,
$$H(m)=-(m+1) \sum_q {1\over q}-(m+1)^2 \sum_q {1\over 2q^2}-
\ldots,\e$$
where the sum over $q$ is as above. 
Now, the second terms and above are convergent, so they can
be replaced by a constant. For the first term, we use eq. (24),
implying that
$$H(m)=-(m+1) \log\log(x)+  y(m)\e$$
where $y(m)$ is some constant.
Exponentiating we find an expression for the sieve product $M(m)$,
$$M(m)=k(m) {x\over \log(x)^{m+1}}\left(1+O({1\over \log(x)})
\right) ,\e$$
where we set the constant $r_m$ to give the asymptotic 
formula eq. (14). For example, we have 
$r_1=0.7931\ldots$ and $r_2=0.7060\ldots$, etc. 

Eq. (32) is an exact result, and it shows that indeed
the PDF's are as conjectured.

The two conjectures eq. (7) and eq. (11) allow us to give a good estimate for the 
occurrence of any multiplet, which is exact, it appears, for large enough numbers.
  
\ack
I wish to thank Dr. Ida Deichaite for pointing my attention to this problem.

\refout
\bye